\documentclass[12pt]{article}

\usepackage{amsmath, epsfig, cite}
\usepackage{amssymb}
\usepackage{amsfonts}
\usepackage{latexsym}
\usepackage{float}
\usepackage{color}

\newtheorem{thm}{Theorem}[section]

\newcommand*{\fl}[2]{\left\lfloor\frac{#1}{#2}\right\rfloor}

\numberwithin{equation}{section}

\newcommand{\qed}{{\hfill$\square$}\medskip}

\begin{document}

\begin{center}
{\Large\bf Congruences on sums of $q$-binomial coefficients}
\end{center}

\vskip 2mm \centerline{Ji-Cai Liu$^1$ and Fedor Petrov$^2$}
\begin{center}
{\footnotesize $^1$Department of Mathematics, Wenzhou University, Wenzhou 325035, PR China\\
{\tt jcliu2016@gmail.com} }
\end{center}

\begin{center}
{\footnotesize $^2$St. Petersburg Department of V. A. Steklov Institute of Mathematics of the Russian Academy of Sciences, St. Petersburg State University, Russian Federation\\
{\tt fedyapetrov@gmail.com} }
\end{center}

\vskip 0.7cm \noindent{\bf Abstract.}
We establish a $q$-analogue of Sun--Zhao's congruence on harmonic sums. Based on this $q$-congruence and a $q$-series identity, we prove a congruence conjecture on sums of central $q$-binomial coefficients, which was recently proposed by Guo. We also deduce a $q$-analogue of a congruence due to Apagodu and Zeilberger from Guo's $q$-congruence.


\vskip 2mm
\noindent{\it MR Subject Classifications}: 11B65, 11A07, 05A10

\section{Introduction}
In combinatorics, the numbers ${2k\choose k}$ are called {\it central binomial coefficients}. In 2011, Sun and Tauraso \cite{st-injt-2011} proved that for any prime $p\ge 5$,
\begin{align}
\sum_{k=0}^{p-1}{2k\choose k}\equiv \left(\frac{p}{3}\right)\pmod{p^2},\label{aa-3}
\end{align}
where $\left(\frac{\cdot}{p}\right)$ denotes the Legendre symbol. Similar congruences on sums of combinatorial sequences have been widely studied (see, for example, \cite{az-amm-2017,st-aam-2010,st-injt-2011}).

In 2017, Apagodu and Zeilberger \cite{az-amm-2017} utilized an algorithm to find and prove some congruences on indefinite sums of many combinatorial sequences such as Catalan numbers, Motzkin numbers and central binomial coefficients. In particular, Apagodu and Zeilberger \cite[Proposition 4]{az-amm-2017} showed that for any prime $p\ge 5$,
\begin{align}
\sum_{i=0}^{p-1}\sum_{j=0}^{p-1}{i+j\choose i}^2\equiv \left(\frac{p}{3}\right)\pmod{p},\label{aa-4}
\end{align}
which they also conjectured holds modulo $p^2$.

In 2010, Sun and Zhao \cite{sz-jnt-2010} proved that for any positive odd integer $n$ and
prime $p\ge n+1$,
\begin{align}
\sum_{0<i_1<\cdots<i_n<p }\left(\frac{i_1}{3}\right)\frac{(-1)^{i_1}}{i_1\cdots i_n}
\equiv 0\pmod{p}. \label{sz-1}
\end{align}
The special case $n=1$ of \eqref{sz-1} reads as
\begin{align}
\sum_{k=1}^{p-1}\left(\frac{k}{3}\right)\frac{(-1)^{k}}{k}
\equiv 0\pmod{p}. \label{sz-2}
\end{align}

$q$-Analogues of congruences (also called {\it $q$-congruences}) were also widely discussed by
several authors (see, for instance, \cite{andrews-dm-1999,guo-ijnt-2018,gs-2018,gw-tjm-2019,gz-aam-2010,gz-am-2019,
Tauraso-aam-2012,tauraso-cm-2013,zudilin-2019}).
Recently, Guo and Zudilin \cite{gz-am-2019} developed a creative microscoping method
to prove many interesting $q$-congruences.

To continue the story of \eqref{aa-3}, \eqref{aa-4} and \eqref{sz-2}, we need some $q$-series notation.
The Gaussian $q$-binomial coefficients are defined as
\begin{align*}
{n\brack k}={n\brack k}_q
=\begin{cases}
\displaystyle\frac{(q;q)_n}{(q;q)_k(q;q)_{n-k}} &\text{if $0\leqslant k\leqslant n$},\\[10pt]
0 &\text{otherwise,}
\end{cases}
\end{align*}
where the $q$-shifted factorial is given by $(a;q)_n=(1-a)(1-aq)\cdots(1-aq^{n-1})$ for $n\ge 1$ and $(a;q)_0=1$. The $q$-integers are defined by $[n]_q=(1-q^n)/(1-q)$.
Moreover, the $n$th cyclotomic polynomial is given by
\begin{align*}
\Phi_n(q)=\prod_{\substack{1\le k \le n\\
(n,k)=1}}(q-e^{2k\pi i/n}).
\end{align*}

In this paper, we first establish a $q$-analogue of \eqref{sz-2} as follows.
\begin{thm}\label{t-1}
For any positive integer $n$, we have
\begin{align}
\sum_{k=1}^{n-1}\left(\frac{k}{3}\right)
\frac{(-1)^k}{[k]_q}q^{\frac{k}{3}\left(k-\left(\frac{k}{3}\right)\right)-\frac{(k-1)(k-2)}{6}}
\equiv \left(\frac{n}{3}\right)\frac{\left(\frac{n}{3}\right)-n}{6}q^{\frac{n^2-1}{3}}(1-q)
\pmod{\Phi_n(q)}.\label{sz-3}
\end{align}
\end{thm}

The congruence \eqref{aa-3} also possesses the following amazing $q$-analogue, which was recently
conjectured by Guo \cite[Conjecture 4.1]{guo-ijnt-2018}.
\begin{thm}\label{t-2}
For any positive integer $n$, we have
\begin{align}
\sum_{k=0}^{n-1}q^k{2k\brack k}\equiv \left(\frac{n}{3}\right)q^{\frac{n^2-1}{3}}\pmod{\Phi_n(q)^2}.\label{a-1}
\end{align}
\end{thm}

Tauraso \cite[Corollary 4.3]{Tauraso-aam-2012} have showed that \eqref{a-1} holds modulo $\Phi_n(q)$. By using \eqref{sz-3}, we shall prove that \eqref{a-1} holds modulo $\Phi_n(q)^2$. Based on \eqref{a-1}, we also give a $q$-analogue of \eqref{aa-4} by establishing the following more general result.
\begin{thm}\label{t-3}
For any positive integer $n$, we have
\begin{align}
\sum_{i=0}^{n-1}\sum_{j=0}^{n-1}q^{j^2+i+j}{i+j\brack i}^2\equiv \left(\frac{n}{3}\right)q^{\frac{n^2-1}{3}}\pmod{\Phi_n(q)^2}.\label{a-2}
\end{align}
\end{thm}

The rest of this paper is organized as follows. Section 2 is devoted to proving Theorem \ref{t-1}. We prove Theorems \ref{t-2} and \ref{t-3} in Sections 3 and 4, respectively.

\section{Proof of Theorem \ref{t-1}}
Note that
\begin{align*}
\sum_{k=1}^{n-1}\left(\frac{k}{3}\right)
\frac{(-1)^k}{[k]_q}q^{\frac{k}{3}\left(k-\left(\frac{k}{3}\right)\right)-\frac{(k-1)(k-2)}{6}}
&=-\sum_{k=0}^{\fl{n-2}{3}}\frac{(-1)^kq^{\frac{3k(k+1)}{2}}}{[3k+1]_q}
-\sum_{k=0}^{\fl{n-3}{3}}\frac{(-1)^kq^{\frac{3k^2+9k+4}{2}}}{[3k+2]_q}\notag\\
&=-\sum_{k=-\fl{n}{3}}^{\fl{n-2}{3}}\frac{(-1)^kq^{\frac{3k(k+1)}{2}}}{[3k+1]_q},
\end{align*}
where $\lfloor x \rfloor$ denotes the integral part of real $x$.
So \eqref{sz-3} is equivalent to
\begin{align}
\sum_{k=-\fl{n}{3}}^{\fl{n-2}{3}}\frac{(-1)^kq^{\frac{3k(k+1)}{2}}}{1-q^{3k+1}}
\equiv \left(\frac{n}{3}\right)\frac{n-\left(\frac{n}{3}\right)}{6}q^{\frac{n^2-1}{3}}\pmod{\Phi_n(q)}.\label{b-1}
\end{align}

We shall distinguish three cases to prove \eqref{b-1}.

{\bf Case 1}\quad $n\equiv 0\pmod{3}$.
This case is equivalent to
\begin{align*}
\sum_{k=-n}^{n-1}\frac{(-1)^kq^{3k(k+1)/2}}{1-q^{3k+1}}\equiv 0\pmod{\Phi_{3n}(q)}.
\end{align*}
Let $\alpha$ be a primitive $3n$th root of unity and
\begin{align*}
C_k=\frac{(-1)^k\alpha^{3k(k+1)/2}}{1-\alpha^{3k+1}}.
\end{align*}
It suffices to show that
\begin{align*}
\sum_{k=-n}^{n-1}C_k=0.
\end{align*}
By $\alpha^{3n(n+1)/2}=(-1)^{n+1}$, we have
\begin{align*}
C_{n-k}=-C_{-k},
\end{align*}
for $1\le k\le n$. It follows that
\begin{align*}
\sum_{k=-n}^{n-1}C_k
=\sum_{k=1}^n\left(C_{n-k}+C_{-k}\right)
=0.
\end{align*}

{\bf Case 2}\quad $n\equiv 1\pmod{3}$. It is equivalent to
\begin{align}
\sum_{k=-n}^{n-1}\frac{(-1)^kq^{3k(k+1)/2}}{1-q^{3k+1}}\equiv
\frac{n}{2}q^{3n^2+2n}\pmod{\Phi_{3n+1}(q)}.\label{b-2}
\end{align}
Letting $k\to k-n-1$ on the left-hand side of \eqref{b-2} and simplifying, we find that \eqref{b-2} is equivalent to
\begin{align}
\sum_{k=1}^{2n}\frac{(-1)^{k}q^{k(3k-1)/2}}{1-q^{3k-1}}\equiv
-\frac{n}{2}\pmod{\Phi_{3n+1}(q)}.\label{b-3}
\end{align}
To prove \eqref{b-3}, it suffices to show that
\begin{align}
\sum_{k=1}^{2n}\frac{(-1)^{k}\zeta^{k(3k-1)/2}}{1-\zeta^{3k-1}}=
-\frac{n}{2},\label{b-4}
\end{align}
where $\zeta$ is a primitive $(3n+1)$th root of unity.

Note that
\begin{align*}
\frac{1}{1-\zeta^{3 k-1}}=\frac12 \left(\frac{1 }{1-\zeta^{(3k-1)/2}}+\frac{1}{1+\zeta^{(3k-1)/2}}\right),
\end{align*}
and
\begin{align*}
\frac{\zeta^{k(3k-1)/2}}{1-\zeta^{(3k-1)/2}}&=\frac{1}{1-\zeta^{(3k-1)/2}}-\sum_{s=0}^{k-1}\zeta^{(3k-1)s/2},\\
\frac{(-1)^k\zeta^{k(3k-1)/2}}{1+\zeta^{(3k-1)/2}}&=\frac{1}{1+\zeta^{(3k-1)/2}}-\sum_{s=0}^{k-1}(-1)^s\zeta^{(3k-1)s/2}.
\end{align*}
Thus,
\begin{align}
\sum_{k=1}^{2n}\frac{(-1)^k\zeta^{k(3k-1)/2}}{1-\zeta^{3k-1}}
&=\frac{1}{2} \sum _{k=1}^{2 n} \left(\frac{(-1)^k }{1-\zeta^{(3k-1)/2}}+\frac{1}{1+\zeta^{(3k-1)/2}}\right)\notag\\
&- \frac{1}{2} \sum _{k=1}^{2 n} \sum_{s=0}^{k-1} \left( (-1)^k+(-1)^s\right)\zeta^{(3k-1)s/2}.
\label{new-1}
\end{align}
Furthermore,
\begin{align}
&\sum _{k=1}^{2 n} \sum_{s=0}^{k-1} \left( (-1)^k+(-1)^s\right)\zeta^{(3k-1)s/2}\notag\\
&=\sum _{s=0}^{2 n-1} \sum_{k=s+1}^{2n} \left( (-1)^k+(-1)^s\right)\zeta^{(3k-1)s/2}\notag\\
&=2n+2\sum_{s=1}^{2n-1}\frac{(-1)^s\zeta^{s(3s+5)/2}}{1-\zeta^{3s}}+\sum_{s=1}^{2n-1}
\left(\frac{\zeta^{(3n+1)s}}{1+\zeta^{3s/2}}-
\frac{(-1)^s\zeta^{(3n+1)s}}{1-\zeta^{3s/2}}\right)\notag\\
&=2n+2\sum_{s=1}^{2n-1}\frac{(-1)^s\zeta^{s(3s+5)/2}}{1-\zeta^{3s}}+\sum_{s=1}^{2n-1}
\left(\frac{1}{1+\zeta^{3s/2}}-
\frac{(-1)^s}{1-\zeta^{3s/2}}\right).\label{new-2}
\end{align}
Letting $s\to 2n+1-k$ on the right-hand side of \eqref{new-2} gives
\begin{align}
&\sum _{k=1}^{2 n} \sum_{s=0}^{k-1} \left( (-1)^k+(-1)^s\right)\zeta^{(3k-1)s/2}\notag\\
&=2n+2\sum_{k=2}^{2n}\frac{(-1)^k\zeta^{k(3k-1)/2}}{1-\zeta^{3k-1}}+\sum_{k=2}^{2n}
\left(\frac{\zeta^{(3k-1)/2}}{1+\zeta^{(3k-1)/2}}-\frac{(-1)^k\zeta^{(3k-1)/2}}{1-\zeta^{(3k-1)/2}}\right)\notag\\
&=2n+2\sum_{k=1}^{2n}\frac{(-1)^k\zeta^{k(3k-1)/2}}{1-\zeta^{3k-1}}+\sum_{k=1}^{2n}
\left(\frac{\zeta^{(3k-1)/2}}{1+\zeta^{(3k-1)/2}}-\frac{(-1)^k\zeta^{(3k-1)/2}}{1-\zeta^{(3k-1)/2}}\right).\label{new-3}
\end{align}
Combining \eqref{new-1} and \eqref{new-3}, we obtain
\begin{align}
\sum_{k=1}^{2n}\frac{(-1)^k\zeta^{k(3k-1)/2}}{1-\zeta^{3k-1}}&=
\frac{1}{4}\sum _{k=1}^{2 n} \left(\frac{(-1)^k(1+\zeta^{(3k-1)/2}) }{1-\zeta^{(3k-1)/2}}+\frac{1-\zeta^{(3k-1)/2}}{1+\zeta^{(3k-1)/2}}\right)-\frac{n}{2}\notag\\
&=-n+\frac{1}{2}\sum_{k=1}^{2n}\left(\frac{(-1)^k}{1-\zeta^{(3k-1)/2}}+\frac{1}{1+\zeta^{(3k-1)/2}}\right).
\label{new-4}
\end{align}

Note that
\begin{align}
&\sum_{k=1}^{2n}\left(\frac{(-1)^k}{1-\zeta^{(3k-1)/2}}+\frac{1}{1+\zeta^{(3k-1)/2}}\right)\notag\\
&=\sum_{k=1}^{2n}\left(\frac{1}{1-\zeta^{(3k-1)/2}}+\frac{1}{1+\zeta^{(3k-1)/2}}\right)
-2\sum_{k=1}^n\frac{1}{1-\zeta^{3k-2}}\notag\\
&=2\sum_{k=1}^{2n}\frac{1}{1-\zeta^{3k-1}}
-2\sum_{k=1}^n\frac{1}{1-\zeta^{3k-2}}.\label{new-5}
\end{align}
Clearly,
\begin{align}
\sum_{k=1}^{2n}\frac{1}{1-\zeta^{3k-1}}
-\sum_{k=1}^n\frac{1}{1-\zeta^{3k-2}}&=\sum_{k=-n+1}^{n}\frac{1}{1-\zeta^{3k-2}}-\sum_{k=1}^n\frac{1}{1-\zeta^{3k-2}}\notag\\
&=\sum_{k=-n+1}^{0}\frac{1}{1-\zeta^{3k-2}}\notag\\
&=\sum_{k=1}^n\frac{1}{1-\zeta^{3k-1}}.\label{new-6}
\end{align}
It follows from \eqref{new-4}--\eqref{new-6} that
\begin{align*}
\sum_{k=1}^{2n}\frac{(-1)^k\zeta^{k(3k-1)/2}}{1-\zeta^{3k-1}}
&=-n+\sum_{k=1}^n\frac{1}{1-\zeta^{3k-1}}\\
&=-n+\frac{1}{2}\sum_{k=1}^n\left(\frac{1}{1-\zeta^{3k-1}}+\frac{1}{1-\zeta^{3(n+1-k)-1}}\right)\\
&=-n+\frac{1}{2}\sum_{k=1}^n\left(\frac{1}{1-\zeta^{3k-1}}+\frac{1}{1-\zeta^{1-3k}}\right)\\
&=-\frac{n}{2}.
\end{align*}

{\bf Case 3}\quad $n\equiv 2\pmod{3}$. Let $\omega$ be a primitive $(3n+2)$th root of unity.
Similarly to the case $n\equiv 1\pmod{3}$, it suffices to show that
\begin{align}
\sum_{k=1}^{2n+1}\frac{(-1)^k\omega^{k(3k+1)/2}}{1-\omega^{3k}}=-\frac{n+1}{2}.\label{b-5}
\end{align}

By using the same method as in the case $n\equiv 1\pmod{3}$, we can get
\begin{align}
\sum_{k=1}^{2n+1}\frac{(-1)^k\omega^{k(3k+1)/2}}{1-\omega^{3k}}
&=\frac{1}{2} \sum _{k=1}^{2 n+1} \left(\frac{(-1)^k\omega^{k/2} }{1-\omega^{3k/2}}+\frac{\omega^{k/2}}{1+\omega^{3k/2}}\right)\notag\\
&- \frac{1}{2} \sum _{k=1}^{2 n+1} \sum_{s=0}^{k-1} \left( (-1)^k+(-1)^s\right)\omega^{(3s+1)k/2}.
\label{b-6}
\end{align}
Clearly,
\begin{align}
&\sum _{k=1}^{2 n+1} \sum_{s=0}^{k-1} \left( (-1)^k+(-1)^s\right)\omega^{(3s+1)k/2}\notag\\
&=\sum _{s=0}^{2 n} \sum_{k=s+1}^{2n+1} \left( (-1)^k+(-1)^s\right)\omega^{(3s+1)k/2}\notag\\
&=2\sum_{s=0}^{2n}\frac{(-1)^s\omega^{(s+2)(3s+1)/2}}{1-\omega^{3s+1}}-\sum_{s=0}^{2n}
\left(\frac{(-1)^s\omega^{(3s+1)(n+1)}}{1-\omega^{(3s+1)/2}}+
\frac{\omega^{(3s+1)(n+1)}}{1+\omega^{(3s+1)/2}}\right).\label{b-7}
\end{align}
Letting $s\to 2n+1-k$ on the right-hand side of \eqref{b-7} gives
\begin{align}
&\sum _{k=1}^{2 n+1} \sum_{s=0}^{k-1} \left( (-1)^k+(-1)^s\right)\omega^{(3s+1)k/2}\notag\\
&=2\sum_{k=1}^{2n+1}\frac{(-1)^k\omega^{k(3k+1)/2}}{1-\omega^{3k}}
-\sum_{k=1}^{2n+1}\left(\frac{(-1)^k\omega^{k/2}}{1-\omega^{3k/2}}
+\frac{\omega^{k/2}}{1+\omega^{3k/2}}\right).\label{b-8}
\end{align}
Combining \eqref{b-6} and \eqref{b-8}, we obtain
\begin{align}
\sum_{k=1}^{2n+1}\frac{(-1)^k\omega^{k(3k+1)/2}}{1-\omega^{3k}}
&=\frac{1}{2} \sum _{k=1}^{2 n+1} \left(\frac{(-1)^k \omega^{k/2} }{1-\omega^{3k/2}}+\frac{\omega^{k/2}}{1+\omega^{3k/2}}\right).\label{b-9}
\end{align}

Since
\begin{align*}
\sum _{k=1}^{2 n+1} \frac{(-1)^k \omega^{k/2} }{1-\omega^{3k/2}}=
-\sum _{k=1}^{2 n+1} \frac{\omega^{k/2} }{1-\omega^{3k/2}}+2\sum _{k=1}^{n} \frac{\omega^{k} }{1-\omega^{3k}},
\end{align*}
we have
\begin{align}
\sum _{k=1}^{2 n+1} \left(\frac{(-1)^k \omega^{k/2} }{1-\omega^{3k/2}}+\frac{\omega^{k/2}}{1+\omega^{3k/2}}\right)
=2\sum _{k=1}^{n} \frac{\omega^{k} }{1-\omega^{3k}}-2\sum_{k=1}^{2n+1}\frac{\omega^{2k}}{1-\omega^{3k}}.\label{b-10}
\end{align}
Now we split the right-hand side of \eqref{b-10} into two identical halfs, one of which reads as
(letting $k\to 3n+2-k$)
\begin{align*}
\sum _{k=1}^{n} \frac{\omega^{k} }{1-\omega^{3k}}-\sum_{k=1}^{2n+1}\frac{\omega^{2k}}{1-\omega^{3k}}
=\sum_{k=n+1}^{3n+1}\frac{\omega^k}{1-\omega^{3k}}-\sum_{k=2n+2}^{3n+1}\frac{\omega^{2k}}{1-\omega^{3k}}.
\end{align*}
Collecting with another half we get
\begin{align*}
&\sum _{k=1}^{2 n+1} \left(\frac{(-1)^k \omega^{k/2} }{1-\omega^{3k/2}}+\frac{\omega^{k/2}}{1+\omega^{3k/2}}\right)\\
&=\sum_{k=1}^{3n+1}\frac{\omega^k-\omega^{2k}}{1-\omega^{3k}}\\
&=\sum_{k=1}^{3n+1}\omega^k\frac{1-\omega^{3k(n+1)}}{1-\omega^{3k}}\\
&=\sum_{k=1}^{3n+1}\left(\omega^k+\omega^{4k}+\omega^{7k}+\cdots+\omega^{(3n+1)k}\right).
\end{align*}
For any integer $d$ with $(3n+2)\nmid d$, we have
\begin{align*}
\sum_{k=0}^{3n+1}\omega^{kd}=0.
\end{align*}
It follows that
\begin{align}
\sum _{k=1}^{2 n+1} \left(\frac{(-1)^k \omega^{k/2} }{1-\omega^{3k/2}}+\frac{\omega^{k/2}}{1+\omega^{3k/2}}\right)=-n-1.\label{b-11}
\end{align}
Then the proof of \eqref{b-5} follows from \eqref{b-9} and \eqref{b-11}.
\qed

\section{Proof of Theorem \ref{t-2}}
We begin with the following identity:
\begin{align}
\sum_{k=0}^{n-1}q^k{2k\brack k}=\sum_{k=0}^{n-1}\left(\frac{n-k}{3}\right)q^{\frac{1}{3}\left(2(n-k)^2-(n-k)
\left(\frac{n-k}{3}\right)-1\right)}{2n\brack k},\label{c-1}
\end{align}
which was conjectured by Z.-W. Sun and proved by Tauraso in a more general form (see \cite[Theorem 4.2]{Tauraso-aam-2012}).
By $1-q^n\equiv 0\pmod{\Phi_n(q)}$, we have
\begin{align*}
1-q^{2n}=(1+q^n)(1-q^n)\equiv 2(1-q^n)\pmod{\Phi_n(q)^2}.
\end{align*}
It follows that for $1\le k\le n-1$,
\begin{align}
{2n\brack k}&=\frac{(1-q^{2n})(1-q^{2n-1})\cdots(1-q^{2n-k+1})}{(1-q)(1-q^2)\cdots(1-q^k)}\notag\\
&\equiv 2(1-q^{n})\frac{(1-q^{-1})\cdots(1-q^{-k+1})}{(1-q)(1-q^2)\cdots(1-q^k)}\pmod{\Phi_n(q)^2}\notag\\
&=2(q^n-1)\frac{(-1)^kq^{-\frac{k(k-1)}{2}}}{1-q^k}.\label{c-2}
\end{align}
Substituting \eqref{c-2} into the right-hand side of \eqref{c-1} gives
\begin{align}
&\sum_{k=0}^{n-1}q^k{2k\brack k}\notag\\
&=\left(\frac{n}{3}\right)q^{\frac{1}{3}\left(2n^2-n
\left(\frac{n}{3}\right)-1\right)}+\sum_{k=1}^{n-1}\left(\frac{n-k}{3}\right)q^{\frac{1}{3}\left(2(n-k)^2-(n-k)
\left(\frac{n-k}{3}\right)-1\right)}{2n\brack k}\notag\\
&\equiv \left(\frac{n}{3}\right)q^{\frac{1}{3}\left(2n^2-n
\left(\frac{n}{3}\right)-1\right)}\notag\\
&+2(q^n-1)\sum_{k=1}^{n-1}(-1)^k\left(\frac{n-k}{3}\right)\frac{q^{\frac{1}{3}\left(2(n-k)^2-(n-k)
\left(\frac{n-k}{3}\right)-1\right)-\frac{k(k-1)}{2}}}{1-q^k}\pmod{\Phi_n(q)^2}.\label{c-3}
\end{align}

Furthermore, letting $k\to n-k$ in the following sum gives
\begin{align}
&\sum_{k=1}^{n-1}(-1)^k\left(\frac{n-k}{3}\right)\frac{q^{\frac{1}{3}\left(2(n-k)^2-(n-k)
\left(\frac{n-k}{3}\right)-1\right)-\frac{k(k-1)}{2}}}{1-q^k}\notag\\
&=\sum_{k=1}^{n-1}(-1)^{n-k}\left(\frac{k}{3}\right)\frac{q^{\frac{1}{3}\left(2k^2-k
\left(\frac{k}{3}\right)-1\right)-\frac{(n-k)(n-k-1)}{2}}}{1-q^{n-k}}.\label{c-4}
\end{align}
For $1\le k\le n-1$, we have
\begin{align}
\frac{q^{\frac{1}{3}\left(2k^2-k
\left(\frac{k}{3}\right)-1\right)-\frac{(n-k)(n-k-1)}{2}}}{1-q^{n-k}}&=
\frac{q^{\frac{1}{3}\left(2k^2-k
\left(\frac{k}{3}\right)-1\right)-\frac{n(n-1)}{2}-\frac{k(k+1)}{2}+nk}}{1-q^{n-k}}\notag\\
&\equiv \frac{(-1)^{n+1}q^{\frac{1}{3}\left(2k^2-k
\left(\frac{k}{3}\right)-1\right)-\frac{k(k+1)}{2}}}{1-q^{-k}}\pmod{\Phi_n(q)}\notag\\
&=\frac{(-1)^{n}q^{\frac{k}{3}\left(k-
\left(\frac{k}{3}\right)\right)-\frac{(k-1)(k-2)}{6}}}{1-q^{k}},\label{c-5}
\end{align}
where we have used the fact that $q^{n(n-1)/2}\equiv (-1)^{n+1}\pmod{\Phi_n(q)}$.

Combining \eqref{c-3}--\eqref{c-5}, we obtain
\begin{align}
\sum_{k=0}^{n-1}q^k{2k\brack k}
&\equiv \left(\frac{n}{3}\right)q^{\frac{1}{3}\left(2n^2-n
\left(\frac{n}{3}\right)-1\right)}\notag\\
&+\frac{2(q^n-1)}{1-q}\sum_{k=1}^{n-1}\left(\frac{k}{3}\right)
\frac{(-1)^k}{[k]_q}q^{\frac{k}{3}\left(k-\left(\frac{k}{3}\right)\right)-\frac{(k-1)(k-2)}{6}}
\pmod{\Phi_n(q)^2}.\label{c-6}
\end{align}
It follows from \eqref{sz-3} and \eqref{c-6} that
\begin{align}
\sum_{k=0}^{n-1}q^k{2k\brack k}\equiv \left(\frac{n}{3}\right)q^{\frac{n^2-1}{3}}\left(q^{\frac{n\left(n-\left(\frac{n}{3}\right)\right)}{3}}
-\frac{n-\left(\frac{n}{3}\right)}{3}(q^n-1)\right)\pmod{\Phi_n(q)^2}.\label{c-8}
\end{align}

By the binomial theorem and $n-\left(\frac{n}{3}\right)\equiv 0\pmod{3}$, we have
\begin{align}
q^{\frac{n\left(n-\left(\frac{n}{3}\right)\right)}{3}}
&=\left(1+q^n-1\right)^{\frac{n-\left(\frac{n}{3}\right)}{3}}\notag\\
&=1+\frac{n-\left(\frac{n}{3}\right)}{3}(q^n-1)+{\frac{n-\left(\frac{n}{3}\right)}{3}\choose 2}(q^n-1)^2+\cdots+(q^n-1)^{\frac{n-\left(\frac{n}{3}\right)}{3}}\notag\\
&\equiv 1+\frac{n-\left(\frac{n}{3}\right)}{3}(q^n-1)\pmod{\Phi_n(q)^2}.\label{c-9}
\end{align}
Then the proof of \eqref{a-1} follows from \eqref{c-8} and \eqref{c-9}.

\section{Proof of Theorem \ref{t-3}}
Let $i+j=m$. Then
\begin{align*}
\sum_{i=0}^{n-1}\sum_{j=0}^{n-1}q^{j^2+i+j}{i+j\brack i}^2
=\sum_{m=0}^{n-1}\sum_{j=0}^{m}q^{j^2+m}{m\brack j}^2
+\sum_{m=n}^{2n-2}\sum_{j=m-n+1}^{n-1}q^{j^2+m}{m\brack j}^2.
\end{align*}
Note that for $n\le m\le 2n-2$ and $m-n+1\le j\le n-1$,
\begin{align*}
{m\brack j}=\frac{(1-q^m)(1-q^{m-1})\cdots(1-q^{m-j+1})}{(1-q)(1-q^2)\cdots(1-q^j)}\equiv 0\pmod{\Phi_n(q)},
\end{align*}
because the numerator contains the factor $1-q^n$ and all of the factors of the denominator are
relatively prime with $\Phi_n(q)$. It follows that
\begin{align*}
\sum_{i=0}^{n-1}\sum_{j=0}^{n-1}q^{j^2+i+j}{i+j\brack i}^2
\equiv\sum_{m=0}^{n-1}\sum_{j=0}^{m}q^{j^2+m}{m\brack j}^2\pmod{\Phi_n(q)^2}.
\end{align*}
Furthermore, by the $q$-Chu-Vandermonde summation formula (see \cite[(3.3.10)]{andrews-b-1998}), we get
\begin{align*}
\sum_{j=0}^{m}q^{j^2}{m\brack j}^2=\sum_{j=0}^{m}q^{j^2}{m\brack j}{m\brack m-j}={2m\brack m}.
\end{align*}
Thus,
\begin{align}
\sum_{i=0}^{n-1}\sum_{j=0}^{n-1}q^{j^2+i+j}{i+j\brack i}^2
\equiv\sum_{m=0}^{n-1}{2m\brack m}\pmod{\Phi_n(q)^2}.\label{d-1}
\end{align}
Then the proof of \eqref{a-2} follows from \eqref{a-1} and \eqref{d-1}.

\vskip 5mm \noindent{\bf Acknowledgments.}
The authors would like to thank Nemo (MathOverflow User) for the helpful discussion and suggestions. The first author was supported by the National Natural Science Foundation of China (grant 11801417).


\begin{thebibliography}{99}
\small \setlength{\itemsep}{-.8mm}

\bibitem{andrews-b-1998}G.E. Andrews, The Theory of Partitions, Cambridge University Press, Cambridge, 1998.

\bibitem{andrews-dm-1999}G.E. Andrews, $q$-Analogs of the binomial coefficient congruences of Babbage, Wolstenholme and Glaisher, Discrete Math. 204 (1999), 15--25.

\bibitem{az-amm-2017}M. Apagodu and D. Zeilberger, Using the ``freshman's dream" to prove combinatorial congruences, Amer. Math. Monthly 124 (2017), 597--608.

\bibitem{guo-ijnt-2018}V.J.W. Guo, Proof of a $q$-congruence conjectured by Tauraso, Int. J. Number Theory (2018), doi: 10.1142/S1793042118501713.

\bibitem{gs-2018}V.J.W. Guo and M.J. Schlosser, Some $q$-supercongruences from transformation formulas for basic hypergeometric series, preprint, 2018, arXiv:1812.06324.

\bibitem{gw-tjm-2019}V.J.W. Guo and S.-D. Wang, Factors of sums and alternating sums of products of $q$-binomial coefficients and powers of $q$-integers, Taiwanese J. Math. 23 (2019), 11--27.

\bibitem{gz-aam-2010}V.J.W. Guo and J. Zeng, Some congruences involving central $q$-binomial coefficients, Adv. Appl. Math. 45 (2010), 303--316.

\bibitem{gz-am-2019}V.J.W. Guo and W. Zudilin, A $q$-microscope for supercongruences, Adv. Math., to appear.

\bibitem{st-aam-2010}Z.-W. Sun and R. Tauraso, New congruences for central binomial coefficients, Adv. Appl. Math. 45 (2010), 125--148.

\bibitem{st-injt-2011}Z.-W. Sun and R. Tauraso, On some new congruences for binomial coefficients, Int. J. Number Theory 7 (2011), 645--662.

\bibitem{sz-jnt-2010}Z.-W. Sun and L.-L. Zhao, Some curious congruences modulo primes,
¡¡ J. Number Theory 130 (2010), 930--935.

\bibitem{Tauraso-aam-2012}R. Tauraso, $q$-Analogs of some congruences involving Catalan numbers, Adv. Appl. Math. 48 (2012), 603--614.
    
\bibitem{tauraso-cm-2013}R. Tauraso, Some $q$-analogs of congruences for central binomial sums,
Colloq. Math. 133 (2013), 133--143.

\bibitem{zudilin-2019}W. Zudilin, Congruences for $q$-binomial coefficients, preprint, 2019, 	arXiv:1901.07843.

\end{thebibliography}
\end{document}